\documentclass[12pt,leqno]{article}
\usepackage{latexsym}
\usepackage{amsmath}
\usepackage{amssymb}

\newcommand{\vf}{\varphi}

\newcommand{\bi}{\begin{itemize}}
\newcommand{\ei}{\end{itemize}}
\newcommand{\beq}{\begin{equation}}
\newcommand{\eeq}{\end{equation}}

\newcommand{\vsp}{\vspace{10mm}}

\newcounter {nrq} 

\begin{document}

\title{\bf\sc Large deviations for random evolutions with independent
 increments in the scheme of L\'{e}vy approximation}

\author{{\sc I.V. Samoilenko}\\
Institute of Mathematics,\\ Ukrainian National Academy of Science,
Kyiv, Ukraine, isamoil@imath.kiev.ua}


\maketitle

{\bf\sc Short title: Large deviations for random evolutions}

\baselineskip 6 mm

\vsp

\hrule
\begin{abstract}
In the work asymptotic analysis of the problem of large deviations
for random evolutions with independent increments in the circuit of
L\'{e}vy approximation is carried out. Large deviations for random
evolutions in the circuit of L\'{e}vy approximation are determined
by exponential generator for jumping process with independent
increments.
\end{abstract}

{\small {\sc Key Words:} {L\'{e}vy approximation, nonlinear
exponential generator, Markov process, locally independent
increments process, piecewise deterministic Markov process, singular
perturbation. }

{\small {\sc Mathematics Subject Classification Primary:} 60J55,
60B10, 60F17, 60K10; Secondary 60G46, 60G60.}

\vsp
\hrule
\section{Introduction}

Asymptotic analysis of the problem of large deviations for random
evolutions with independent increments in the circuit of L\'{e}vy
approximation (see Koroliuk and Limnios, 2005, Ch. 9) is carried out
in the paper.

Asymptotic analysis of random evolutions with independent increments
in L\'{e}vy approximation scheme is conducted in the work of
Koroliuk, Limnios and Samoilenko (2009).

In the monograph of Feng and Kurtz (2006) an effective method for
studying the problem of large deviations for Markov processes is
developed. It is based on the theory of convergence of exponential
(nonlinear) operators. The exponential operator in the series scheme
with a small series parameter $ \varepsilon \to0 (\varepsilon> 0) $
has the form (see, e.g. Koroliuk, 2011):
$$
\mathbb {H} ^ \varepsilon \varphi (x): = e ^ {- \varphi (x) /
\varepsilon} \varepsilon \mathbb {L} ^ \varepsilon e ^ {\varphi (x)
/ \varepsilon},
$$
where the operators $ \mathbb {L} ^ \varepsilon, \varepsilon> 0 $
define Markov processes $ \zeta^ \varepsilon (t), t \geq0,
\varepsilon> 0 $ in the series scheme on the standard phase-space
$(G, \mathcal{G})$. Test-functions $\varphi(x)\in G$ are real-valued
and finite.

Random evolutions with independent increments (see Koroliuk and
Limnios, 2005, Ch. 1) are given by:
$$ \xi (t) = \xi_0 + \int ^ t_0 \eta (ds; x (s)), \ t\geq0. \eqno (1)
$$

Markov processes with independent increments $\eta(t; x), $ $
t\geq0, $ $ x\in E,$ are defined in $\mathbb{R}$ and given by the
generators

$$
\Gamma(x)\varphi(u)=\int_{\mathbb{R}}[\varphi(u+v)-\varphi(u)]\Gamma(dv;
x), \ x\in E, \varphi(u)\in\mathcal{B}_{\mathbb{R}}.
$$

\textbf{Remark 1.} The process in $\mathbb{R}^d, d>1$ may also be
studied. See Remark 7 for more details.

Markov switching process $ x(t), t \geq0, $ on a standard
 phase space $ (E, \mathcal{E}) $ is defined by the generator
$$
Q\varphi(x)=q(x)\int_{\mathbb{E}}[\varphi(y)-\varphi(x)]P(x, dy), \
x\in E, \varphi(u)\in\mathcal{B}_{E}. \eqno(2)
$$

Thus, random evolution (1) is characterized by the generator of
two-component Markov process $\xi(t), x(t), t\geq 0$ (see Koroliuk
and Limnios, 2005, Ch. 2)
$$
\mathbb{L}\varphi(u, x)=Q\varphi(\cdot, x)+\Gamma(x)\varphi(u,
\cdot).
$$

The basic assumption about the switching Markov process is the
following condition

\begin{itemize}
\item[{\bf C1}:]
   Markov process  $x(t), t\geq0,$ is uniformly ergodic with the stationary distribution $\pi(A), \ A\in \mathcal{E}.$
 \end{itemize}

Let $\Pi$ be a projector onto null-subspace of reducible-invertible
operator $Q$, defined in (2):
$$\Pi\varphi(x)=\int_E\pi(dx)\varphi(x).$$

The following correlation is true $$Q\Pi=\Pi Q=0.$$

Potential operator $R_0$ has the following property (Koroliuk and
Limnios, 2005, Ch. 1):
$$QR_0=R_0Q=\Pi-I.$$

\textbf{Remark 2.} It follows from the last correlation that under
solvability condition
$$\Pi\psi=0$$ Poisson
equation
$$Q\varphi=\psi$$ has the unique solution $$\psi=R_0\varphi,$$ when
$\Pi\varphi=0.$

\textbf{Remark 3.} Studying of limit properties of Markov processes
is based at the martingale characterization of such processes,
namely we should regard
$$
\mu_t=\varphi(x(t))-\varphi(x(0))-\int^t_0\mathbb{L}\varphi(x(s))ds,
\eqno(3)
$$
where $\mathbb{L}$ is the generator that defines Markov process
$x(t), t\geq0,$ on the standard phase-space $(E, \mathcal{E})$. It
has a dense domain $\mathcal{D}(\mathbb{L})\subseteq\mathcal{B}_E$,
that contains continuous functions with continuous derivatives. Here
$\mathcal{B}_E$ - Banach space of real-valued finite test-functions
$\varphi(x)\in E,$ endowed by the norm: $\|\varphi\|:=\sup_{x\in
E}|\varphi(x)|.$

Large deviation theory is based on the studying of exponential
martingale characterization (see Feng and Kurtz, 2006, Ch.1):
$$
\widetilde{\mu}_t=\exp\{\varphi(x(t))-\varphi(x(0))-\int^t_0\mathbb{H}\varphi(x(s))ds\}
 \eqno(4)
$$ is the martingale.

Here exponential nonlinear operator
$$
\mathbb{H}\varphi(x):=e^{-\varphi(x)}\mathbb{L}e^{\varphi(x)}, \
\varphi(x)\in \mathcal{B}_E.
$$

Equivalence of (3) and (4) follows from the correlations:

\textit{Proposition} (see Ethier and Kurtz, 1986, p.66)
$$\mu(t)=x(t)-\int_0^ty(s)ds$$ is the martingale if and only if
$$\widetilde{\mu}(t)=x(t)exp\left\{-\int_0^t\frac{y(s)}{x(s)}ds\right\} \mbox{is the martingale.}$$
We may assume that the domain $\mathcal{D}(\mathbb{L})$ contains
constants, and if $\varphi(x)\in\mathcal{D}(\mathbb{L})$, then there
exists a constant $c$ such that
$\varphi(x)+c\in\mathcal{D}(\mathbb{L})$ is positive.

\textbf{Remark 4.} The large deviation problem is realized in four
stages (Feng and Kurtz, 2006, Ch.2):

1) Verify the convergence of the exponential (nonlinear) generator
that defines large deviations;

2) Verify the exponential tightness of Markov processes;

3) Verify the comparison principle for the limit exponential
generator;

4) Construct a variational representation for the limit exponential
generator.

The stages 2)--4) for the exponential generator corresponding to the
processes with independent increments are realized in Feng and Kurtz
(2006). Some of the stages are also presented in the monograph of
Freidlin and Wentzel (1998), where the large deviation problem is
studied with the use of cumulant of the process with independent
increments. Cumulant and exponential generator are obviously
connected.

Really, generator of Markov process may be written in the form (see,
e.g. Skorokhod, 1989)
$$\mathbb{L}\varphi(x)=\int_\mathbb{R}e^{\lambda
x}a(\lambda)\overline{\varphi}(\lambda)d\lambda,$$ where
$a(\lambda)$ - cumulant of the process,
$\overline{\varphi}(\lambda)=\int_\mathbb{R}e^{-\lambda
x}\varphi(x)dx.$

Inverse transformation gives
$$\int_\mathbb{R}e^{-\lambda
x}\mathbb{L}\varphi(x)dx=a(\lambda)\overline{\varphi}(\lambda).$$

Let's rewrite $$\int_\mathbb{R}e^{-\lambda
x}\mathbb{L}\varphi(x)dx=\int_\mathbb{R}e^{-\lambda
x}a(\lambda)\varphi(x)dx,$$ and by taking $$e^{-\lambda
x}\varphi(x)=:\widetilde{\varphi}(x)$$ we obtain
$$\int_\mathbb{R}e^{-\lambda x}\mathbb{L}e^{\lambda
x}\widetilde{\varphi}(x)dx=\int_\mathbb{R}a(\lambda)\widetilde{\varphi}(x)dx.$$

Thus, $$e^{-\lambda x}\mathbb{L}e^{\lambda x}=a(\lambda),$$ or,
using the exponential generator:
$$\mathbb{H} \varphi_0(x)=a(\lambda), \hskip2mm \mbox{where} \hskip2mm \varphi_0(x)=\lambda x.$$

Our aim is to realize the stage 1) - to verify the convergence of
the exponential (nonlinear) generator that defines large deviations
$\mathbb{H}^{\varepsilon,\delta}\varphi^\delta_\varepsilon(u,x):=e^{-\varphi^\delta_\varepsilon/\varepsilon}\varepsilon
\mathbb{L}^{\delta}_{\varepsilon}
e^{\varphi^\delta_\varepsilon/\varepsilon}$ (see Theorem 1):
$$\mathbb{H}^{\varepsilon,\delta}\varphi^\delta_\varepsilon(u,x)\to H^0\varphi(u), \varepsilon, \delta\to 0, \varepsilon^{-1}\delta\to 1.$$ To do this we use the method of solution of the problem of
singular perturbation with two small series parameters.

Normalization of random evolution (1) by a small series parameters
for solution of large deviation problem in L\'{e}vy approximation
scheme is realized in a following way:
$$\xi_{\varepsilon}^{\delta}(t)=\xi_{\varepsilon}^{\delta}(0)+\int^t_0\eta^{\delta}_{\varepsilon}(ds; x(s/\varepsilon^3)), \ t\geq0,$$
$$\eta^{\delta}_{\varepsilon}(t)=\varepsilon\eta^{\delta}(t/\varepsilon^3),$$
$$\Gamma^{\delta}_{\varepsilon}(x)\varphi(u)=\varepsilon^{-3}\int_{\mathbb{R}}[\varphi(u+\varepsilon v)-\varphi(u)]\Gamma^{\delta}(dv; x), \ x\in
E,$$ where $\varepsilon,\delta\to 0$ so that
$\varepsilon^{-1}\delta\to 1$.

\textbf{Remark 5.} In the paper (Koroliuk, 2011) V.S.Koroliuk
proposed to use the method of solution of the problem of singular
perturbation for the studying of large deviations for random
evolutions with independent increments in asymptotically small
diffusion scheme.

In classical works asymptotical analysis of the problem of large
deviations is made, as a rule, with the use of large series
parameter $n\to\infty,$ sometimes even some different parameters
(see, e.g. Mogulskii, 1993).

The method, proposed in this work, with the use of two small
parameters, was firstly realized in Samoilenko (2011) for the scheme
of Poisson approximation.

\section {L\'{e}vy approximation conditions}

\begin{itemize}

\item[{\bf C2:}] {\it L\'{e}vy approximation}. The family of processes with independent increments
$\eta^{\delta}(t;x), $ $x\in E,$ $ t\geq 0 $ satisfies L\'{e}vy
approxiation conditions:

\begin{description}
\item[LA1] Approximation of mean values:
$$a_{\delta}(x) = \int_{\mathbb{R}} v\Gamma^{\delta}(dv; x)
= \delta a_1(x)+\delta^2[a(x) +\theta_a^{\delta} (x)],$$ and
$$c_{\delta}(x) = \int_{\mathbb{R}}
v^2\Gamma^{\delta}(dv; x) = \delta^2[c(x) + \theta_c^{\delta}
(x)],$$ where $$\sup\limits_{x\in E}|a_1(x)|\leq
a_1<+\infty,\sup\limits_{x\in E}|a(x)|\leq a<+\infty,
\sup\limits_{x\in E}|c(x)|\leq c<+\infty.$$

\item[LA2] Asymptotic representation of intensity kernel
$$\Gamma_g^{\delta}(x) = \int_{\mathbb{R}} g(v)\Gamma^{\delta}(dv; x)
= \delta^2[\Gamma_g(x) + \theta^{\delta}_g(x)]$$ for all $g \in
C_3(\mathbb{R})$ - measure-determining class of functions (see Jacod
and Shiryaev, 1987, Ch. 7), $\Gamma_g(x)$ is a finite kernel
$$|\Gamma_g(x)| \leq\Gamma_g \quad \hbox{(constant depending on $g$)}.$$

Kernel $\Gamma^0(dv; x)$ is defined on the measure-determining class
of functions $C_3(\mathbb{R})$ by a relation
$$\Gamma_g(x) = \int_{\mathbb{R}} g(v)\Gamma^0(dv; x),\quad g \in
C_3(\mathbb{R}).$$

Negligible terms $\theta_a^\delta,\theta_c^\delta, \theta_g^\delta$
satisfy the condition
$$\sup\limits_{x\in E} |\theta_{\cdot}^{\delta}(x)|\to 0,\quad
\delta\to 0.$$

\item[LA3] Balance condition:
$$\int_E \pi(dx)a_1(x) = 0.$$

\end{description}

\item[{\bf C3}:] {\it Uniform square integrability}:
$$\lim\limits_{c\to\infty}\sup\limits_{x\in E} \int_{|v|>c} v^2\Gamma^0(dv; x) = 0.$$

\item[{\bf C4}:] {\it Exponential finiteness}: $$\int_{\mathbb{R}}e^{p |v|}\Gamma^{\delta}(dv; x)<\infty, \forall p\in \mathbb{R}.$$

\end{itemize}

\section {Main result}

\textbf{ Theorem 1.} Solution of large deviation problem for random
evolution
$$\xi^{\delta}_{\varepsilon}(t)=\xi^{\delta}_{\varepsilon}(0)+\int^t_0\eta^{\delta}_{\varepsilon}(ds;
x(s/\varepsilon^3)), \ t\geq0,$$ defined by a generator of
two-component Markov process $\xi(t), x(t), t\geq 0$ $$
\mathbb{L}^{\delta}_{\varepsilon}\varphi(u,
x)=\varepsilon^{-3}Q\varphi(\cdot,
x)+\Gamma^{\delta}_{\varepsilon}(x)\varphi(u, \cdot), \eqno(5)
$$ where
$$\Gamma^{\delta}_{\varepsilon}(x)\varphi(u)=\varepsilon^{-3}\int_{\mathbb{R}}[\varphi(u+\varepsilon v)-\varphi(u)]\Gamma^{\delta}(dv; x), \ x\in
E \eqno(6)$$ is realized by the exponential generator
$$
H^0\varphi(u)=(\widetilde{a}-\widetilde{a}_0)\varphi'(u)+\frac{1}{2}\sigma^2(\varphi'(u))^2+\int_{\mathbb{R}}[e^{v\varphi'(u)}-1]\widetilde{\Gamma}^0(dv),
\eqno(7)
$$
$$
\widetilde{a}=\Pi a(x)=\int_E\pi(dx)a(x), \widetilde{a}_0=\Pi
a_0(x)=\int_E\pi(dx)a_0(x),
a_0(x)=\int_{\mathbb{R}}v\Gamma^0(dv;x),$$ $$ \widetilde{c}=\Pi
c(x)=\int_E\pi(dx)c(x), \widetilde{c}_0=\Pi
c_0(x)=\int_E\pi(dx)c_0(x),
c_0(x)=\int_{\mathbb{R}}v^2\Gamma^0(dv;x),$$
$$
\sigma^2=(\widetilde{c}-\widetilde{c}_0)+2\int_E\pi(dx)
a_1(x)R_0a_1(x), \widetilde{\Gamma}^0(v)=\Pi
{\Gamma}^0(v;x)=\int_E\pi(dx){\Gamma}^0(v;x).$$

\textbf{Remark 6.} Large deviations for random evolutions in
L\'{e}vy approximation scheme are determined by exponential
generator for jumping process with independent increments. Studying
of large deviation problem for jumping process with independent
increments is presented in monograph (Freidlin and Wentzel, 1998,
Ch.3,4).

\textbf{Remark 7.} The limit exponential generator in the Euclidean
space $\mathbb{R}^d, d>1$ is represented in the following view:
$$
H^0\varphi(u)=\sum^d_{k=1}(\widetilde{a}_k-\widetilde{a}_k^0)\varphi'_k+\frac{1}{2}\sum^d_{k,r=1}\sigma_{kr}\varphi'_k\varphi'_r+\int_{\mathbb{R}^d}[e^{v\varphi'(u)}-1]\widetilde{\Gamma}^0(dv),
\ \varphi'_k:=\partial\varphi(u)/\partial u_k, 1\leq k\leq d.
$$
Here $\sigma^2=[\sigma_{kr}; 1\leq k,r\leq d]$ is the variance
matrix.

In addition, the last exponential generator can be extended on the
space of absolutely continuous functions (see Feng and Kurtz, 2006)
$$C^1_b(R^d)=\{\vf: \ \exists \lim_{|u|\to\infty}\vf(u)=\vf(\infty),
\ \lim_{|u|\to\infty}\vf'(u)=0\}.$$

\textbf{Proof.} Limit transition in the exponential nonlinear
generator of random evolution is realized on the perturbed
test-functions
$$ \varphi^\delta_\varepsilon(u,
x)=\varphi(u)+\varepsilon\ln[1+\delta\varphi_1(u,
x)+\delta^2\varphi_2(u, x)],
$$ where $\varphi(u)\in C^3(\mathbb{R})$ (the space of continuous bounded functions with continuous bounded derivatives up to third degree). Thus, we have from (5): $$
\mathbb{H}^{\varepsilon,\delta}\varphi^\delta_\varepsilon=e^{-\varphi^\delta_\varepsilon/\varepsilon}\varepsilon
\mathbb{L}^{\delta}_{\varepsilon}
e^{\varphi^\delta_\varepsilon/\varepsilon}=e^{-\varphi^\delta_\varepsilon/\varepsilon}
[\varepsilon^{-2}Q+\varepsilon\Gamma_{\varepsilon}^\delta(x)]
e^{\varphi^\delta_\varepsilon/\varepsilon}=$$
$$e^{-\varphi/\varepsilon}[1 +
\delta\varphi_1+\delta^2\varphi_2]^{-1}[\varepsilon^{-2}Q+\varepsilon\Gamma_{\varepsilon}^\delta(x)]e^{\varphi/\varepsilon}[1
+ \delta\varphi_1+\delta^2\varphi_2].
$$ To obtain the asymptotic behavior of the last exponential
generator we use the following results.

\textbf{Lemma 1.} Exponential generator
$$H^{\varepsilon}_{Q}\varphi_\varepsilon^\delta(u,x)=e^{-\varphi_\varepsilon^\delta/\varepsilon}\varepsilon^{-2}Qe^{\varphi_\varepsilon^\delta/\varepsilon}\eqno(8)$$ has the following asymptotic
representation
$$
H^{\varepsilon}_{Q}\varphi_\varepsilon^\delta=\varepsilon^{-1}Q\varphi_1+Q\varphi_2-\varphi_1Q\varphi_1
+\theta_{Q}^{\varepsilon,\delta}(x), \eqno(9)
$$ where $\sup\limits_{x\in
E}|\theta_{Q}^{\varepsilon,\delta}(x)|\to 0, \varepsilon,\delta\to
0.$

\textbf{ Proof.} We have:
$$
\mathbb{H}^{\varepsilon}_Q\varphi^\delta_\varepsilon=e^{-\varphi/\varepsilon}[1
+
\delta\varphi_1+\delta^2\varphi_2]^{-1}\varepsilon^{-2}Qe^{\varphi/\varepsilon}[1
+ \delta\varphi_1+\delta^2\varphi_2]=$$
$$\left[1-\delta\varphi_1+\delta^2\frac{\varphi_1^2+\delta\varphi_1\varphi_2-\varphi_2}{1+\delta\varphi_1+\delta^2\varphi_2}\right]
[\delta
\varepsilon^{-2}Q\varphi_1+\delta^2\varepsilon^{-2}Q\varphi_2]=\delta
\varepsilon^{-2}Q\varphi_1+\delta^2\varepsilon^{-2}Q\varphi_2-\delta^2\varepsilon^{-2}\varphi_1Q\varphi_1+\theta_{Q}^{\varepsilon,\delta}(x),
$$ where $$\theta_{Q}^{\varepsilon,\delta}(x)=\delta^3\varepsilon^{-2}\frac{\varphi_1^2+\delta\varphi_1\varphi_2-\varphi_2}{1+\delta\varphi_1+\delta^2\varphi_2}[Q\varphi_1+\delta Q\varphi_2]-\delta^3\varepsilon^{-2}\varphi_1Q\varphi_2.$$
By the limit condition $\varepsilon^{-1}\delta\to 1,
\varepsilon,\delta\to 0$, we finally have (9).

Lemma is proved.

\textbf{Lemma 2.} Exponential generator
$$H^{\varepsilon,\delta}_{\Gamma}(x)\varphi_\varepsilon^\delta(u,x)=e^{-\varphi_\varepsilon^\delta/\varepsilon}\varepsilon\Gamma^{\delta}_{\varepsilon}(x)e^{\varphi_\varepsilon^\delta/\varepsilon}\eqno(10)$$  has the following
asymptotic representation
$$
H^{\varepsilon,\delta}_{\Gamma}(x)\varphi_\varepsilon^\delta=H_{\Gamma}(x)\varphi(u)+\varepsilon^{-1}
a_1(x)\varphi'(u)+\theta_{\Gamma}^{\varepsilon,\delta}(x),
$$ where $$H_{\Gamma}(x)\varphi(u)=(a(x)-a_0(x))\varphi'(u)+\frac{1}{2}(c(x)-c_0(x))(\varphi'(u))^2+\int_{\mathbb{R}}[e^{v\varphi'(u)}-1]\Gamma^0(dv;x),\eqno(11)$$ and $\sup\limits_{x\in
E}|\theta_{\Gamma}^{\varepsilon,\delta}(x)|\to 0,
\varepsilon,\delta\to 0.$

\textbf{ Proof.} We have:
$$
\mathbb{H}^{\varepsilon,\delta}_{\Gamma}(x)\varphi^\delta_\varepsilon=e^{-\varphi/\varepsilon}[1
+
\delta\varphi_1+\delta^2\varphi_2]^{-1}\varepsilon\Gamma_\varepsilon^\delta(x)e^{\varphi/\varepsilon}[1
+ \delta\varphi_1+\delta^2\varphi_2]=$$
$$e^{-\varphi/\varepsilon}\left[1-\delta\varphi_1+\delta^2\frac{\varphi_1^2+\delta\varphi_1\varphi_2-\varphi_2}{1+\delta\varphi_1+\delta^2\varphi_2}\right]
[\varepsilon\Gamma_\varepsilon^\delta(x)e^{\varphi/\varepsilon}+\varepsilon\delta\Gamma_\varepsilon^\delta(x)e^{\varphi/\varepsilon}\varphi_1+
\varepsilon\delta^2\Gamma_\varepsilon^\delta(x)e^{\varphi/\varepsilon}\varphi_2]=
$$ $$H^{\varepsilon,\delta}_{\Gamma}(x)\varphi(u)+e^{-\varphi/\varepsilon}\varepsilon\delta[\Gamma_\varepsilon^\delta(x)e^{\varphi/\varepsilon}\varphi_1-
\varphi_1\Gamma_\varepsilon^\delta(x)e^{\varphi/\varepsilon}]+\widetilde{\theta}_{\Gamma}^{\varepsilon,\delta}(x),
$$ where $$\widetilde{\theta}_{\Gamma}^{\varepsilon,\delta}(x)=\varepsilon\delta^2[e^{-\varphi/\varepsilon}\Gamma_\varepsilon^\delta(x)e^{\varphi/\varepsilon}\varphi_2-e^{-\varphi/\varepsilon}\varphi_1\Gamma_\varepsilon^\delta(x)e^{\varphi/\varepsilon}\varphi_1]
+\varepsilon\delta^2\frac{\varphi_1^2+\delta\varphi_1\varphi_2-\varphi_2}{1+\delta\varphi_1+\delta^2\varphi_2}
[e^{-\varphi/\varepsilon}\Gamma_\varepsilon^\delta(x)e^{\varphi/\varepsilon}+$$
$$e^{-\varphi/\varepsilon}\delta\Gamma_\varepsilon^\delta(x)e^{\varphi/\varepsilon}\varphi_1+
e^{-\varphi/\varepsilon}\delta^2\Gamma_\varepsilon^\delta(x)e^{\varphi/\varepsilon}\varphi_2]-
\varepsilon\delta^3e^{-\varphi/\varepsilon}\varphi_1\Gamma_\varepsilon^\delta(x)e^{\varphi/\varepsilon}\varphi_2$$

We use the following results:

\textbf{Lemma 3.}
$$\Gamma_{\varepsilon}^{\delta}(x)e^{\varphi(u)/\varepsilon}\varphi_1(u,x)=\varphi_1(u,x)\Gamma_{\varepsilon}^{\delta}(x)e^{\varphi(u)/\varepsilon}+(\varepsilon\delta)^{-1}\widehat{\theta}_{\Gamma}^{\varepsilon,\delta}(x),$$
where for the negligible term $\sup\limits_{x\in
E}|\widehat{\theta}^{\varepsilon,\delta}_\Gamma(x)|\to 0,
\varepsilon,\delta\to 0.$

\textbf{Proof.} Really, by (6) we have:
$$\Gamma_{\varepsilon}^{\delta}(x)e^{\varphi(u)/\varepsilon}\varphi_1(u,x)=\varepsilon^{-3}
\int_{\mathbb{R}}[e^{\varphi(u+\varepsilon
v)/\varepsilon}\varphi_1(u+\varepsilon
v,x)-e^{\varphi(u)/\varepsilon}\varphi_1(u,x)]\Gamma^\delta(dv;x)=$$
$$\varphi_1(u,x)\Gamma_{\varepsilon}^{\delta}(x)e^{\varphi(u)/\varepsilon}+
(\varepsilon\delta)^{-1}\left[\varphi'_1(u,x)\varepsilon^{-1}\delta
\int_{\mathbb{R}}e^{\varphi(u+\varepsilon
v)/\varepsilon}v\Gamma^\delta(dv;x)\right].$$

Let's estimate the last integral. As soon as function $\varphi(u)$
is bounded, we have for fixed $\varepsilon$:
$$\int_{\mathbb{R}}e^{\varphi(u+\varepsilon
v)/\varepsilon}v\Gamma^\delta(dv;x)<e^C\int_{\mathbb{R}}v\Gamma^\delta(dv;x)=\delta
e^C[a_1(x)+\delta a(x)+\delta\theta_a^{\delta}(x)].$$

Thus, we see that the last term is negligible when $\varepsilon,
\delta\to 0$.

Lemma is proved.

\textbf{Lemma 4.} Exponential generator
$$H^{\varepsilon,\delta}_{\Gamma}(x)\varphi(u)=e^{-\varphi/\varepsilon}\varepsilon\Gamma_\varepsilon^\delta(x)e^{\varphi/\varepsilon}\eqno(12)$$
has the following asymptotic representation
$$H^{\varepsilon,\delta}_{\Gamma}(x)\varphi^\delta_\varepsilon=H_\Gamma(x)\varphi(u)+
\varepsilon^{-1}a_1(x)\varphi'(u)+\theta^{\varepsilon,\delta}(x),$$
where $\sup\limits_{x\in E}|\theta^{\varepsilon,\delta}(x)|\to 0,
\varepsilon,\delta\to 0.$

\textbf{Proof.} Let's rewrite (12), using the view of the generator
(6). We have:
$$H^{\varepsilon,\delta}_{\Gamma}(x)\varphi(u)=\varepsilon^{-2}\int_{\mathbb{R}}[e^{\Delta_{\varepsilon}\varphi(u)}-1]\Gamma^\delta(dv;x),
$$ where $$
\Delta_\varepsilon\varphi(u):=\varepsilon^{-1}[\varphi(u+\varepsilon
v)-\varphi(u)].
$$

We may rewrite it in a following way:
$$H^{\varepsilon,\delta}_{\Gamma}(x)\varphi(u)=\varepsilon^{-2}\int_{\mathbb{R}}[e^{\Delta_{\varepsilon}\varphi(u)}-1-\Delta_{\varepsilon}\varphi(u)-
\frac{1}{2}(\Delta_{\varepsilon}\varphi(u))^2]\Gamma^\delta(dv;x)+$$
$$\varepsilon^{-2}\int_{\mathbb{R}}[\Delta_{\varepsilon}\varphi(u)+
\frac{1}{2}(\Delta_{\varepsilon}\varphi(u))^2]\Gamma^\delta(dv;x).
$$

Easy to see that the function
$\psi^{\varepsilon}_u(v)=e^{\Delta_{\varepsilon}\varphi(u)}-1-\Delta_{\varepsilon}\varphi(u)-
\frac{1}{2}(\Delta_{\varepsilon}\varphi(u))^2$ belongs to the class
$C_3(\mathbb{R})$. Really,
$$\psi^{\varepsilon}_u(v)/v^2\to 0, v\to 0.$$ Besides, this function is continuous and bounded for every $\varepsilon$
 under the condition that $\varphi(u)$ is bounded. Moreover, the function
$\psi^{\varepsilon}_u(v)$ is bounded uniformly by $u$ under the
conditions {\bf C3, C4} and if $\varphi'(u)$ is bounded.

Thus, we have:
$$H^{\varepsilon,\delta}_{\Gamma}(x)\varphi(u)=\varepsilon^{-2}\delta^2\int_{\mathbb{R}}[e^{\Delta_{\varepsilon}\varphi(u)}-1-\Delta_{\varepsilon}\varphi(u)-
\frac{1}{2}(\Delta_{\varepsilon}\varphi(u))^2]\Gamma^0(dv;x)+$$
$$\varepsilon^{-2}\int_{\mathbb{R}}[\Delta_{\varepsilon}\varphi(u)-v\varphi'(u)-
\varepsilon
\frac{v^2}{2}\varphi''(u)]\Gamma^\delta(dv;x)+\varepsilon^{-2}\delta
a_1(x)\varphi'(u)+\varepsilon^{-2}\delta^2 a(x)\varphi'(u)+$$
$$\varepsilon^{-1}\delta^2
c(x)\varphi''(u)+
\varepsilon^{-2}\int_{\mathbb{R}}[\frac{1}{2}(\Delta_{\varepsilon}\varphi(u))^2-\frac{v^2}{2}(\varphi'(u))^2]\Gamma^\delta(dv;x)+\varepsilon^{-2}\delta^2\frac{1}{2}
c(x)(\varphi'(u))^2.$$

The functions in the second and third integrals are obviously belong
to $C_3(\mathbb{R})$. Using Taylor's formula to the test-functions
$\varphi(u)\in C^3(\mathbb{R})$, and condition \textbf{PA2} we
obtain:
$$H^{\varepsilon,\delta}_{\Gamma}(x)\varphi(u)=\varepsilon^{-2}\delta^2\int_{\mathbb{R}}[e^{v\varphi'(u)}-1-v\varphi'(u)-
\frac{v^2}{2}(\varphi'(u))^2]\Gamma^0(dv;x)+$$
$$\varepsilon^{-2}\delta^2\int_{\mathbb{R}}(e^{v\varphi'(u)}\varepsilon\frac{v^2}{2}\varphi''(\widetilde{u})-\varepsilon\frac{v^2}{2}\varphi''(\widetilde{u})-
\varepsilon^2\frac{v^4}{8}(\varphi''(\widetilde{u}))^2)\Gamma^0(dv;x)+$$
$$\varepsilon^{-2}\delta^2\int_{\mathbb{R}}\varepsilon^{2}\frac{v^3}{3!}\varphi'''(\widetilde{u})\Gamma^0(dv;x)+\varepsilon^{-2}\delta
a_1(x)\varphi'(u)+\varepsilon^{-2}\delta^2
a(x)\varphi'(u)+\varepsilon^{-1}\delta^2 c(x)\varphi''(u)+$$ $$
\varepsilon^{-2}\delta^2\int_{\mathbb{R}}\varepsilon^{2}\frac{v^4}{4}(\varphi''(\widetilde{u}))^2\Gamma^0(dv;x)+\varepsilon^{-2}\delta^2\frac{1}{2}
c(x)(\varphi'(u))^2.$$

By the limit condition $\varepsilon^{-1}\delta\to 1$,  we finally
have:
$$
H^{\varepsilon,\delta}_{\Gamma}(x)\varphi(u)=H_{\Gamma}(x)\varphi(u)+\varepsilon^{-1}
a_1(x)\varphi'(u)+\theta^{\varepsilon,\delta}(x),$$ where
$\sup\limits_{x\in E}|\theta^{\varepsilon,\delta}(x)|\to 0,
\varepsilon,\delta\to 0.$

Lemma is proved.

From Lemmas 3,4 we obtain
$$
H^{\varepsilon,\delta}_{\Gamma}(x)\varphi_\varepsilon^\delta=H_{\Gamma}(x)\varphi(u)+\varepsilon^{-1}
a_1(x)\varphi'(u)+\theta_{\Gamma}^{\varepsilon,\delta}(x),$$ where
$\sup\limits_{x\in E}|\theta_{\Gamma}^{\varepsilon,\delta}(x)|\to 0,
\varepsilon,\delta\to 0.$

Lemma is proved.

From (8) and (10) we see that
$$\mathbb{H}^{\varepsilon,\delta}\varphi^\delta_\varepsilon=H^{\varepsilon}_{Q}\varphi_\varepsilon^\delta(u,x)+H^{\varepsilon,\delta}_{\Gamma}(x)\varphi_\varepsilon^\delta(u,x)$$

Thus, using Lemmas 1,2, we obtain asymptotic representation $$
\mathbb{H}^{\varepsilon,\delta}\varphi^\delta_\varepsilon=\varepsilon^{-1}[
Q\varphi_1+a_1(x)\varphi'(u)]+Q\varphi_2-\varphi_1Q\varphi_1+H_{\Gamma}(x)\varphi(u)+
h^{\varepsilon,\delta}(x),
$$ where $h^{\varepsilon,\delta}(x)=\theta^{\varepsilon,\delta}_Q(x)+\theta^{\varepsilon,\delta}_{\Gamma}(x).$

Now we may use the solution of singular perturbation problem for
reducibly-invertible operator $Q$ (see Koroliuk and Limnios, 2005,
Ch. 1).
$$\begin{array}{c}
    Q\varphi_1+a_1(x)\varphi'(u)=0, \\
    Q\varphi_2-\varphi_1Q\varphi_1+H_{\Gamma}(x)\varphi(u)=
    H^0\varphi(u).
  \end{array}
 $$

From the first equation we obtain
$$\varphi_1(u,x)=R_0a_1(x)\varphi'(u),\hskip 5mm Q\varphi_1(u,x)=-a_1(x)\varphi'(u).$$ After substitution to the
second equation we have
$$Q\varphi_2+a_1(x)R_0a_1(x)(\varphi'(u))^2+H_{\Gamma}(x)\varphi(u)=
    H^0\varphi(u)$$ From the solvability condition:
 $$H^0\varphi(u)=\Pi
H_{\Gamma}(x)\Pi\varphi(u)+\Pi
a_1(x)R_0a_1(x)\mathbf{1}(\varphi'(u))^2,$$ where $\mathbf{1}$ is
the unit vector.

Now, using (11) we finally obtain (7).

The negligible term $h^{\varepsilon,\delta}(x)$ may be found
explicitly, using the solution of Poisson equation (see Remark 1, in
details Koroliuk and Limnios, 2005)
$$\varphi_2(u,x)=R_0\widetilde{H}(x)\varphi(u)-R_0a_1(x)R_0a_1(x)\mathbf{1}(\varphi'(u))^2, \hskip 5mm
\widetilde{H}(x):=H^0-H_{\Gamma}(x).$$

Theorem is proved.

\end{document}